\documentclass[english,11pt,dvips]{article}

\usepackage{amsmath}
\usepackage{amsfonts}
\usepackage{amssymb}
\usepackage{longtable}

\newtheorem{theorem}{Theorem}

{\bf}{\rm}

\def\Real{\mathbb{R}}
\def\Co{\mathbb{C}}
\def\g{\mathfrak{g}}
\def\h{\mathfrak{h}}

\def\so{\mathfrak{so}}
\def\spin{\mathfrak{spin}}

\def\sl{\mathfrak{sl}}
\def\gl{\mathfrak{gl}}
\def\su{\mathfrak{su}}
\def\u{\mathfrak{u}}
\def\sp{\mathfrak{sp}}

\def\R{\mathcal{R}}
\def\RIC{\mathcal{R}ic}
\def\P{\mathcal{P}}
\def\Q{\mathcal{Q}}
\def\id{\mathop\text{\rm id}\nolimits}
\def\Id{\mathop\text{\rm Id}\nolimits}

\def\spa{\mathop\text{{\rm span}}\nolimits}
\def\Hom{\mathop\text{\rm Hom}\nolimits}

\def\tr{\mathop\text{\rm tr}\nolimits}

\def\Ric{\mathop{{\rm Ric}}\nolimits}
\def\grad{\mathop{{\rm grad}}}
\def\tRic{\mathop{\widetilde{\rm Ric}}\nolimits}

\def\owedge{\wedge}
\def\wJ{\wedge_J}
\def\owJ{\wedge_J}

\newcommand{\be}{\begin{equation}}
\newcommand{\ee}{\end{equation}}

\begin{document}

\title{Covariant derivative of the curvature tensor of  pseudo-K\"ahlerian manifolds}


\author{Anton S. Galaev}

\maketitle

\begin{abstract}

It is well known that the curvature tensor of a pseudo-Riemannian
manifold can be decomposed   with respect to the pseudo-orthogonal
group  into the sum of the Weyl conformal curvature tensor, the traceless part of the
Ricci tensor and of the scalar  curvature. A similar decomposition
with respect to the pseudo-unitary group  exists on a
pseudo-K\"ahlerian manifold; instead of the Weyl tensor one
obtains the Bochner tensor.  In the present paper, the known
decomposition with respect to the pseudo-orthogonal group of the
covariant derivative of the curvature tensor of a
pseudo-Riemannian manifold is refined. A decomposition with
respect to the pseudo-unitary group of the covariant derivative of
the curvature tensor for pseudo-K\"ahlerian manifolds is obtained.
This defines natural classes of spaces generalizing locally symmetric spaces
and Einstein spaces. It is shown that the values of the covariant
derivative of the curvature tensor for a non-locally symmetric
pseudo-Riemannian manifold with an irreducible connected holonomy
group different from the pseudo-orthogonal and pseudo-unitary
groups belong to an irreducible module of the holonomy group.



\end{abstract}

\section{Introduction}

The well-known Ricci decomposition of the curvature tensor $R$ of
a pseudo-Riemannian manifold $(M,g)$ defines the following three
components of $R$:
 \begin{description}
 \item the Weyl conformal curvature tensor $W$, which is the totally trace-free  part of
$R$;
\item the trace-free part $\Ric^0$ of the Ricci tensor $\Ric$ of $(M,g)$;
\item the scalar curvature $s$.
\end{description}

This decomposition is the consequence of the fact that the space
$\R(\so(p,q))$ of algebraic curvature tensors (i.e. the space of
possible values of the curvature tensor of a pseudo-Riemannian
manifold $(M,g)$) decomposes into the direct sum of three
irreducible $\so(p,q)$-modules if $p+q\geq 5$, where $(p,q)$ is
the signature of $g$.
 Setting to zero some
of these components gives rise to different geometrical
characterizations of~$(M,g)$:

\begin{center}
\begin{tabular}{rcl} $W=0$& $\Leftrightarrow$ & $(M,g)$ is conformally flat ($p+q\geq 4$);\\
 $\Ric^0=0$ & $\Leftrightarrow$ &   $(M,g)$ is an Einstein manifold ($p+q\geq 3$);\\
$\Ric^0=0$ and $s=0$ & $\Leftrightarrow$ &   $(M,g)$ is Ricci-flat;\\
$W=0$ and $\Ric^0=0$ & $\Leftrightarrow$ &   $(M,g)$ has constant
section curvature ($p+q\geq 3$).
\end{tabular}\end{center}

In \cite{Gray} Gray showed that the space $\RIC^\nabla$ of
possible values of the covariant derivative of the Ricci tensor of
$(M,g)$ decomposes into the direct sum of three irreducible
$\so(p,q)$-modules. Setting to zero some of the corresponding
components of the  covariant derivative of the Ricci tensor gives
6 natural systems of equations, that are generalizations of the
Einstein equation. Chapter 16 from \cite{Besse} is dedicated to
these equations.

In \cite{Str88} Strichartz decomposed the space
$\R^\nabla(\so(p,q))$ of covariant derivatives of algebraic
curvature tensors (i.e. the space of possible values of the
covariant derivative of the curvature tensor of a
pseudo-Riemannian manifold $(M,g)$) into four components. He used
the $\so(p,q)$-equivariant linear map
$$\tr_{2,4}:\R^\nabla(\so(p,q)) \to \RIC^\nabla,$$
which has the form $$\nabla R\mapsto \nabla\Ric.$$ The kernel of
this  map is an irreducible submodule of $\R^\nabla(\so(p,q))$ and
$\tr_{2,4}$ restricted to the other three irreducible submodules
of $\R^\nabla(\so(p,q))$ is an isomorphism. Then an inverse
isomorphism to this one is constructed.

As the first result of this paper, we refine the decomposition
from \cite{Str88}. For this we consider the $\so(p,q)$-equivariant
linear map
$$\tr_{1,5}:\R^\nabla(\so(p,q)) \to \P(\so(p,q)).$$
Here the space $\P(\h)$ for a subalgebra $\h\subset\so(p,q)$ is
defined in the following way:
$$\P(\h)=\{P\in\Hom(\Real^{p,q},\h)|g(P(X)Y,Z)+g(P(Y)Z,X)+g(P(Z)X,Y)=0,\quad X,Y,Z\in\Real^{p,q}\},$$
where $g$ is the pseudo-Euclidean form on $\Real^{p,q}$. The
results of \cite{onecomp} allow to find spaces $\P(\h)$ for each
irreducible holonomy algebra $\h\subset\so(p,q)$ of
pseudo-Riemannian manifolds. In particular, the $\so(p,q)$-module
$\P(\so(p,q))$ is the direct sum of two irreducible modules. The
kernel of $\tr_{1,5}$  consists of two irreducible submodules of
$\R^\nabla(\so(p,q))$ and $\tr_{1,5}$ restricted to the other two
irreducible submodules of $\R^\nabla(\so(p,q))$ is an isomorphism.
We construct the  inverse isomorphism to this one and find the
projections of $\nabla R$ to the irreducible parts of the kernel
of $\tr_{1,5}$.

We show that the covariant derivative $\nabla R$ of the curvature
tensor $R$ of a pseudo-Riemannian manifold can be decomposed into
the following four components:

\begin{description} \item
the totally trace-free  part of $\nabla R$, which coincides with
the totally trace-free  part of $\nabla W$\\ (we denote this
component by $S'_0$); \item the symmetrization  of the tensor
$\nabla\left(\Ric-\frac{2}{n+2}sg\right)$
 (we denote this component by $S''_0$); \item
the Cotton curvature tensor $C$; \item the gradient $\grad s$ of
the scalar curvature $s$. \end{description}

Comparing with the decomposition from \cite{Str88}, the obtained
below components have a simpler form. Especially, we found  the
explicit form for $S'_0$, while in \cite{Str88} it is written only
that $S'_0$ equals $\nabla R$ minus the other three components.
The fact that $S'_0$ belongs to the space $\R^\nabla(\so(p,q))$
implies the second Bianchi identity for $\nabla W$ \cite{Eis97}.

Remark that in \cite{Str88,TV81} an unpublished preprint
\cite{Gray-V} of Gray and Vanhecke is cited, where the
decomposition of the  space $\R^\nabla(\so(p,q))$ is also found.

The decomposition of the space $\R^\nabla(\so(1,n))$ was applied
in \cite{2sym} to the classification of Lorentzian spaces
satisfying the condition $\nabla^2 R=0$.

Now one may set to zero some of the obtained  components of $\nabla
R$ and this will give 14 different natural systems of  equations.
The last three components  of $\nabla R$ are defined exactly by
the same tensors defining the three components of $\nabla \Ric$
and may be obtained taking a certain trace in the decomposition of
$\nabla R$. Hence 6 of the obtained 14 systems of equations are considered
in Chapter 16 from \cite{Besse}.

Recall the geometric meaning of some of these equations:

\begin{tabular}{rcl} $C=0$& $\Leftrightarrow$ & $W$ is harmonic;\\
 $C=0$ and $\grad s=0$ & $\Leftrightarrow$ & $R$ is harmonic $\Leftrightarrow$ $W$
 is harmonic and $s$ is constant\\ & & $\Leftrightarrow$ $\Ric$ is a Codazzi tensor;\\
$\nabla_X\Ric(X,X)=0$ & $\Leftrightarrow$ &   $S''_0=0$ and $\grad s=0$;\\
$\nabla_X\left(\Ric(X,X)-\frac{2}{n+2}sg(X,X)\right)=0$ &
$\Leftrightarrow$ & $S_0''=0$.
\end{tabular}

In \cite{DR77,DR09} and in some other papers the equation $$\nabla
W=0$$ is studied. In particular, if $(M,g)$ is a Riemannian
manifold, then this equation implies $W=0$ or $\nabla R=0$
\cite{DR77}. The final step in a local classification of
pseudo-Riemannian manifolds with $\nabla W=0$ is done in
\cite{DR09}. In terms of the decomposition of $\nabla R$ the
equation $\nabla W=0$ is equivalent to the following two equations:
$$S_0'=0,\qquad C=0.$$ It would be interesting to describe
pseudo-Riemannian manifolds $(M,g)$ that satisfy the stronger
condition $$S_0'=0.$$

Then we consider the case of pseudo-K\"ahlerian manifolds. It is
well known that in this case the  curvature tensor $R$  can be
decomposed into the following three components:
 \begin{description}
 \item the Bochner tensor $B$, which is the totally trace-free  part of
$R$;
\item the trace-free part $\Ric^0$ of the Ricci tensor $\Ric$ of $(M,g)$;
\item the scalar curvature $s$.
\end{description}

We introduce some notation simplifying the usual formulas for the
decomposition of $R$.

Pseudo-K\"ahlerian manifolds with $B=0$ are called
Bochner-K\"ahler. Many results and references concerning these
manifolds can be found in the fundamental paper of Bryant
\cite{Bryant2001}.

Then we show that the covariant derivative $\nabla R$ of the
curvature tensor of a  pseudo-K\"ahlerian manifold can be
decomposed into the following three components:

\begin{description} \item
the totally trace-free  part of $\nabla R$, which coincides with
the totally trace-free  part of $\nabla B$\\ (we denote this
component by $Q_0$);
 \item a tensor $D$, which is an analog of
the Cotton curvature tensor;
\item the gradient $\grad s$ of the scalar curvature $s$.
\end{description}

To obtain this decomposition we use our technics developed for
general pseudo-Riemannian manifolds.

This decomposition allows to consider 6 natural systems of
equations on $\nabla R$. The pseudo-K\"ahlerian manifolds with
$D=0$ (i.e. with harmonic Bochner tensor) are studied e.g. in
\cite{Kim09}. The pseudo-K\"ahlerian manifolds with $\nabla B=0$
(i.e. with $Q_0=D=0$) are studied e.g. in \cite{Mat69,Mat73}.

The fact that $Q_0$ belongs to the space $\R^\nabla(\u(p,q))$
implies the second Bianchi identity for $\nabla B$, which is
obtained recently in \cite{Omachi03}.

In \cite{Omachi91} the tensor $\nabla R$ for  Riemannian and
K\"ahlerian manifolds  is decomposed in two different ways in the
sum of two orthogonal components, this gives some inequalities for
the norms of $\nabla R$, $\nabla\Ric$ and $\grad s$. The results
of the present paper applied to Riemannian and K\"ahlerian
manifolds show that    $|\nabla R|^2$ equals to the squares of the
norms of the components of $\nabla R$. This implies different
inequalities (that can be obtained omitting the norms of some
components), some of these inequalities  are already obtained in
\cite{Omachi91}. If some of the inequalities become an equality on
a manifold $(M,g)$, then the omitted components of $\nabla R$ must
be zero. 

In Sections \ref{secR}, \ref{secnR}, \ref{secRK} and \ref{secnRK}
we recall some facts and state the results. In Sections
\ref{secPTh1}, \ref{secPTh2}, \ref{secPTh3} and \ref{secPTh4} we
prove the results. Finally in Section \ref{secOther} we show that
the values of the covariant derivative of the curvature tensor for
a non-locally symmetric pseudo-Riemannian manifold with an
irreducible holonomy algebras different from $\so(p,q)$ and
$\u(p,q)$ belong to an irreducible module of the holonomy algebra.

There are other approaches to algebraic covariant derivatives of
curvature tensors, see e.g. \cite{DFGG} and the references
therein.

Below the letters $X,Y,Z,U,V$ denote either elements of
$\Real^{p,q}$, or of $\Real^{2p,2q}$, or vector fields on $(M,g)$.
Similarly, $X_1,...,X_n$ is either a basis of $\Real^{p,q}$ or a
local basis of vector fields on $(M,g)$. In the pseudo-K\"ahlerian
case $X_1,...,X_n,X_{n+1}=JX_1,...,X_{2n}=JX_{n}$ is either a
basis of $\Real^{2p,2q}$ or a local basis of vector fields on
$(M,g,J)$.

{\bf Acknowledgement.} The author is grateful to
Dmitri~V.~Alekseevsky and Rod Gover for useful discussions.

\section{Decomposition of the curvature tensor of a pseudo-Riemannian manifold}\label{secR}

Denote be $g$ the pseudo-Euclidean metric on $\Real^{p,q}$. Let
$n=p+q$. Using $g$ we may identify the Lie algebra $\so(p,q)$ with
the space of bivectors $\Lambda^2\Real^{p,q}$ in such a way that
$$(X\wedge Y)Z= g (X,Z)Y- g (Y,Z)X.$$  For an element
$H\in\odot^2\Real^{p,q}$ (which can be considered both as a
symmetric linear map and as a symmetric bilinear form) denote by
$H\owedge g$ the endomorphism of $\wedge^2\Real^{p,q}$ defined by
\begin{equation} (H\owedge g)(X,Y)=HX\wedge Y+X\wedge
HY.\end{equation} The same notation can be used for a
pseudo-Riemannian manifold $(M,g)$.

The curvature tensor $R$ of any pseudo-Riemannian manifold $(M,g)$
can be decomposed into the  sum of three components in the
following way:
\begin{equation}\label{decR} R=W-\frac{1}{n-2}\left(\Ric-\frac{s}{n}g\right)\owedge g-\frac{s}{2n(n-1)}g\owedge g,
\end{equation}
where $\Ric$ and $s$ are the  Ricci tensor and the scalar
curvature, respectively. The tensor $W$ is called the Weyl
conformal curvature tensor.

Decomposition \eqref{decR} can be rewritten as
\begin{equation}\label{decR1} R=W-L\owedge g,\end{equation}
where $$L=\frac{1}{n-2}\left(\Ric-\frac{s}{2(n-1)}g\right)$$ is
the Schouten tensor.

Let us explain the origin of Decomposition \eqref{decR}. Let
$(p,q)$ be the signature of the manifold $(M,g)$, $p+q=n$. For any
subalgebra $\h\subset\so(p,q)$ consider the space of algebraic
curvature tensors
$$\R(\h)=\{R\in \Hom(\Lambda^2\Real^{p,q},\h)|R(X,Y)Z+R(Y,Z)X+R(Z,X)Y=0\}.$$
In \cite{Al} the spaces $\R(\h)$ are found for all irreducible
Riemannian holonomy algebras $\h\subset\so(n)$. This result allows
to find also the spaces $\R(\h)$ for all irreducible holonomy
algebras $\h\subset\so(p,q)$ of pseudo-Riemannian manifolds. In
particular, if $n\geq 5$, then  the $\so(p,q)$-module
$\R(\so(p,q))$ admits the following decomposition into irreducible
components:
 \begin{equation}\label{decRso} \R(\so(p,q))=\R_0(\so(p,q))\oplus\R'(\so(p,q))\oplus\R_1(\so(p,q))\simeq V_{2\bar\pi_2}\oplus V_{2\pi_1}\oplus \Real,\end{equation}
where  $\bar\pi_2=2\pi_2$ if $n=5$ and $\bar\pi_2=\pi_2$ if $n\geq
6$; here and in the next section $V_\Lambda$ denotes the real
irreducible representation of $\so(p,q)$ defined by the complex
irreducible representation of $\so(n,\Co)$ with the highest weight
$\Lambda$. For simple Lie algebras we use the notations from
\cite{V-O}. If $x\in M$ and $R$ is the curvature tensor of
$(M,g)$, then, clearly, $R_x\in \R(\so(T_xM))\simeq\R(\so(p,q))$.
Decomposition \eqref{decRso} shows that $R_x$ can be written as
the sum
$$R_x=R_{0x}+R'_x+R_{1x}$$ where $R_{0x}$, $R'_x$ and $R_{1x}$ belong to the
submodules of
 $\R(\so(T_x,M))$ isomorphic
 to $V_{2\bar\pi_2}$,  $V_{2\pi_1}$ and $\Real$, respectively.
 Note that $V_{2\pi_1}\oplus
\Real\simeq \odot^2\Real^{p,q}$ and the inclusion
$\odot^2\Real^{p,q}\hookrightarrow \R(\so(p,q))$ is given by
$$H\in \odot^2\Real^{p,q}\quad\mapsto\quad H\owedge g
\in\R(\so(p,q)).$$ In particular, $\Real\subset\R(\so(p,q))$
coincides with $\Real g\owedge g$. Decomposition \eqref{decR} at
the point $x\in M$ has the form
$$R_x=(R_x-H_{0x}\owedge g_x-\mu g_x\owedge g_x)+H_{0x}\owedge g_x+\mu g_x\owedge
g_x,$$ where $\mu\in\Real$ and $H_{0x}$ belongs to the submodule
of $\odot^2 T_xM$ isomorphic to $V_{2\pi_1}$, i.e. $H_{0x}$ is
trace-free. Note that the map
$$\R(\so(p,q))\to\odot^2\Real^{p,q},\quad R\mapsto \Ric(R),\quad \Ric(R)(X,Y)=\tr(Z\mapsto R(Z,X)Y)$$
is $\so(p,q)$-equivariant, hence it is zero on the submodule
$V_{2\bar\pi_2}\subset\R(\so(p,q))$. Similarly, if $R\in
V_{2\pi_1}$, then $g^{ij}\Ric(X_i,X_j)=0$. These conditions
 imply $H_{0x}=-\frac{1}{n-2}\left(\Ric_x-\frac{s_x}{n}g_x\right)$ and $\mu=-\frac{s_x}{2n(n-1)}$.
 Omitting the point $x$ in the above decomposition
of $R_x$, we get Decomposition \eqref{decR}.

If $n=3$, then
\begin{equation}\label{decRson3} \R(\so(p,q))\simeq  V_{4\pi_1}\oplus \Real,\end{equation}
consequently, $W=0$. If $(p,q)\in\{(4,0),(2,2),(0,4)\}$, then
\begin{equation}\label{decRson4} \R(\so(p,q))\simeq (V_{4\pi_1}\oplus V_{4\pi'_1})\oplus V_{2\pi_1+2\pi'_1}\oplus \Real,\end{equation}
and $W$ admits a decomposition
\begin{equation}\label{Wn4} W=W^++W^-\end{equation} into the
self-dual and anti self-dual parts. Recall that the Cotton tensor
$C$ is defined as follows:
$$C(X,Y,Z)=(n-2)\big(\nabla_ZL(X,Y)-\nabla_YL(X,Z)\big).$$
It holds $$(n-3)C(X,Y,Z)=(n-2)g^{ij}\nabla_{X_i}W(Y,Z,X,X_j).$$ If
$(p,q)\in\{(4,0),(2,2),(0,4)\}$, then this and  \eqref{Wn4} define
the decomposition
\begin{equation}\label{Cn4} C=C^++C^-.\end{equation}

\section{Decomposition of the covariant derivative of the curvature
tensor of a pseudo-Riemannian manifold}\label{secnR}

For any subalgebra $\h\subset\so(p,q)$ define the space
$$\R^\nabla(\h)=\{S\in\Hom\big(\Real^{p,q},\R(\h)\big)|S_X(Y,Z)+S_Y(Z,X)+S_Z(X,Y)=0\}.$$
From the second Bianchi identity it follows that if $(M,g)$ is a
pseudo-Riemannian manifold with the holonomy algebra $\h$ at a
point $x\in M$, and $\nabla R$ is the covariant derivative of the
curvature tensor of $(M,g)$, then $\nabla R_x\in\R^\nabla(\h)$.

First we find the space $\R^\nabla(\so(p,q))$.

\begin{theorem}\label{Th1} The $\so(p,q)$-module
  $\R^\nabla(\so(p,q))$ admits the following decomposition into the sum of
 irreducible $\so(p,q)$-modules:
\begin{align} \R^\nabla(\so(p,q))&\simeq V_{6\pi_1}\oplus V_{4\pi_1}\oplus \Real^{p,q},   \quad\text{ if } n=3,\\
\R^\nabla(\so(p,q))&\simeq    (V_{5\pi_1+\pi'_1}\oplus
V_{\pi_1+5\pi'_1})\oplus
V_{3\pi_1+3\pi'_1}\oplus(V_{3\pi_1+\pi'_1}\oplus
V_{\pi_1+3\pi'_1})\oplus\Real^{p,q},\\ \nonumber &  \quad\text{ if } (p,q)\in\{(4,0),(2,2),(0,4)\},\\
\R^\nabla(\so(p,q))&\simeq V_{\pi_1+2\bar\pi_2}\oplus
V_{3\pi_1}\oplus V_{\pi_1+\bar\pi_2}\oplus\Real^{p,q},\quad\text{
if } n\geq 5,
\end{align} where $\bar\pi_2=2\pi_2$ if $n=5$ and $\bar\pi_2=\pi_2$ if $n\geq 6$.
\end{theorem}

Now we give the explicit form of the above decomposition for the
covariant derivative $\nabla R$  of the  curvature tensor $R$ of a
pseudo-Riemannian manifold  $(M,g)$ of signature $(p,q)$. First we
set some notation. Define the following tensors:
\begin{align}g(H_XY,Z)=&\frac{1}{3(n+1)}\Big(C(Y,Z,X)+C(Z,Y,X)\Big),\\
g(T_{0X}Y,Z)=&-\frac{1}{3(n-2)}\Big(\nabla_X\Ric(Y,Z)+\nabla_Y\Ric(Z,X)+\nabla_Z\Ric(X,Y)
\\ \nonumber &
-\frac{2}{n+2}\big(g(\grad s, X)g(Y,Z) +g(\grad s, Y)g(Z,X)+
g(\grad s, Z)g(X,Y)\big)\Big),\\
g(T_{1X}Y,Z)=&-\frac{1}{2(n-1)(n+2)}\big(g(\grad s,
X)g(Y,Z)+g(\grad s, Y)g(Z,X)+
g(\grad s, Z)g(X,Y)\big),\\
g(P_0(X)Y,Z)=&\frac{1}{3(n+1)}C(X,Y,Z),\\
(\varphi_1)_X=&H_X\owedge g,\\
(\varphi_2)_X(Y,Z)=&P_0\big((Y\wedge Z)X\big)+X\wedge
\big(P_0(Z)Y-P_0(Y)Z\big).\end{align} In particular, $T_0$ and
$T_1$ are the symmetrizations of the tensors
$-\frac{1}{(n-2)}\left(\nabla\Ric-\frac{2}{n+2}sg\right)$ and\\
$-\frac{3}{2(n-1)(n+2)}g(\grad s,\cdot)g$, respectively. Note that
$T=T_0+T_1$, where
$$g(T_XY,Z)=-\frac{1}{3}\left(\nabla_XL(Y,Z)+\nabla_YL(Z,X)+\nabla_ZL(X,Y)\right)$$
is the symmetrization of the tensor $-\nabla L$.

\begin{theorem}\label{Th2} Let $(M,g)$ be a pseudo-Riemannian manifold of signature $(p,q)$, $p+q=n$. Then the covariant derivative $\nabla R$
of the  curvature tensor $R$ of $(M,g)$ admits the following
decomposition:
\begin{equation}\label{decnablaR} \nabla R=S_0'+S_0''+S'+S_1,\end{equation}
where  \begin{align}
\label{S0'}S_0'&=\nabla W+\frac{3}{n-2}\varphi_1-3\varphi_2,\\
(S_0'')_X&=T_{0X}\owedge g,\\
\label{S'}S'&=\varphi_1+3\varphi_2,\\
(S_1)_X&=T_{1X}\owedge g.\end{align} If $n=3$, then $S_0'=0$ and
$\varphi_1=\varphi_2$. If $(p,q)\in\{(4,0),(2,2),(0,4)\}$, then
$S_0'$ and $S'$ can be further decomposed:
$$S_0'={S_0'}^++{S_0'}^-,\quad  S'={S'}^++{S'}^-,$$ where ${S_0'}^\pm$ and ${S'}^\pm$ are given
 by the same formulas as $S'_0$ and $S'$, respectively, with $W$ and
$C$ replaced by $W^\pm$ and $C^{\pm}$.
\end{theorem}

The fact that $S_0'$ satisfies the second Bianchi identity implies
the known second Bianchi identity for the tensor $W$ \cite{Eis97},
we may rewrite it in the form \begin{multline}\label{2BW}
g\big((\nabla_X W(Y,Z)+\nabla_Y W(Z,X)+\nabla_Z
W(X,Y))V,U\big)\\=-\frac{1}{n-2}\Big(C(U,X,Y)g(Z,V)+C(U,Y,Z)g(X,V)+C(U,Z,X)g(Y,V)\\-C(V,X,Y)g(Z,U)
-C(V,Y,Z)g(X,U)-C(V,Z,X)g(Y,U)\Big).\end{multline}

Note that $$(S_0''+S_1)_X=T_X\owedge g.$$

Consider some traces for the obtained tensors. The tensor $S'_0$ is totally trace-free. It holds
\begin{align}  (S''_0)_{X_i}(Z,X_j)g^{ij}&=0,\\
 (S')_{X_i}(X,Y,Z,X_j)g^{ij}&=C(Z,X,Y),\\
 (S_1)_{X_i}(X,Y,Z,X_j)g^{ij}&=\frac{1}{2(n-1)}\big(g(\grad
s,X)g(Y,Z)-g(\grad s,Y)g(X,Z)\big),\\
 (S''_0)_{X}(X_i,Y,X_j,Z)g^{ij}&=-(n-2)g(T_{0X}Y,Z),\\
(S')_{X}(X_i,Y,X_j,Z)g^{ij}&=\frac{1}{3}\big(C(Y,Z,X)+C(Z,Y,X)\big),\\
(S_1)_{X}(X_i,Y,X_j,Z)g^{ij}&=\frac{1}{2(n-1)}g(\grad
s,X)g(Y,Z)-(n-2)g(T_{1X}Y,Z).
\end{align}
In \cite{Gray,Besse} it is shown that the space $\RIC^\nabla_x$ of
possible values of the tensor $\nabla\Ric_x$ admits a
decomposition
$$\RIC^\nabla_x=Q_x\oplus S_x\oplus A_x.$$
into the sum of irreducible $\so(T_xM)$-modules. This corresponds
to the decomposition
$$\nabla\Ric=\xi_Q+\xi_S+\xi_A.$$
It can be checked that $$\xi_Q=\frac{1}{2(n-1)}g(\grad s,\cdot)
g-(n-2)T_1,\quad \xi_S=-(n-2)T_0,\quad
(\xi_A)_X(Y,Z)=\frac{1}{3}\big(C(Y,Z,X)+C(Z,Y,X)\big).$$ This
shows that the decomposition of $\nabla\Ric$ from
\cite{Gray,Besse} can be obtained using the above decomposition of
$\nabla R$ taking $\tr_{(2,4)}$.

\section{Decomposition of  the curvature
tensor of a pseudo-K\"ahlerian manifold}\label{secRK}

Let $(M,g,J)$ be a pseudo-K\"ahlerian manifold of signature
$(2p,2q)$. Let $n=p+q\geq 2$. In order to simplify the usual
expression for the decomposition of the curvature tensor of
$(M,g,J)$ we set some notation, similar ideas can be found in
\cite{ACG}. The tangent space to $(M,g,J)$ is identified with the
pseudo-Euclidean space $\Real^{2p,2q}$ endowed with a
pseudo-Euclidean metric $g$ and a $g$-orthogonal complex structure
$J$. Using $J$, we may identify $\Real^{2p,2q}$ with $\Co^n$ in
such a way that $J$ corresponds to the multiplication by the
complex unit $\mathrm{i}$. Consider the corresponding
pseudo-Hermitian metric
$$h(X,Y)=g(X,Y)+g(X,JY)\mathrm{i},\quad X,Y\in\Real^{2p,2q}.$$
The expression \be (X\wJ Y)Z=h(Z,X)Y-h(Z,Y)X\ee defines an element
$X\wJ Y\in\u(p,q)$. Note that $$X\wJ Y=X\wedge Y+JX\wedge JY.$$
This construction allows to identify $\u(p,q)$ with the space
$$\wedge^2_J\Real^{2p,2q}=\spa\{X\wJ
Y|X,Y\in\Real^{2p,2q}\}\subset\wedge^2\Real^{2p,2q}.$$ Under this
identification, the complex structure $J\in\u(p,q)$ corresponds to
the element $$\frac{1}{2}\sum_{i=1}^n\epsilon_ie_i\wJ Je_i,$$
where $e_1,...,e_n,Je_1,...,Je_n$ is an orthonormal basis (i.e.
 it is an orthogonal basis such that $g(e_i,e_i)=\epsilon_i$, where
$\epsilon_i=-1$ for $i=1,...,p$ and $\epsilon_i=1$ for
$i=p+1,...,n$). The complex space $\wedge^2\Co^n$ may be
considered as the following subspace of $\wedge^2\Real^{2p,2q}$:
$$\wedge^2\Co^n=\{X\wedge Y-JX\wedge JY|X,Y\in\Real^{2p,2q}\}.$$
We obtain the orthogonal decomposition
$$\wedge^2\Real^{2p,2q}=\wedge_J^2\Real^{2p,2q}\oplus\wedge^2\Co^n=\u(p,q)\oplus\wedge^2\Co^n.$$
Let $R\in\R(\u(p,q))$. Then
$$R:\wedge^2\Real^{2p,2q}\to\u(p,q)\subset\wedge^2\Real^{2p,2q}$$
is a symmetric linear map. Consequently, $R$ is zero on the
orthogonal complement to $\u(p,q)$ in $\wedge^2\Real^{2p,2q}$,
i.e. $R$ is zero on $\wedge^2\Co^n.$ This implies the known
equalities
\begin{equation} R(JX,JY)=R(X,Y),\quad
R(JX,Y)+R(X,JY)=0.\end{equation} Note that \be R(X\wJ Y)=R(X\wedge
Y)+R(JX\wedge JY)=2R(X\wedge Y).\ee Note also that $$R(X\wedge
Y)=R(X\otimes Y-Y\otimes X)=2R(X,Y).$$ Thus, \be R(X\wJ
Y)=4R(X,Y).\ee

From \cite{Al} it follows that the space $\R(\u(p,q))$ admits the
following decomposition:
$$\R(\u(p,q))=\R_0(\u(p,q))\oplus\R'(\u(p,q))\oplus\R_1(\u(p,q)),$$
where $\R_0(\u(p,q))\otimes\Co\simeq V_{2\pi_1+2\pi_{n-1}}$ and
$$\R'(\u(p,q))\oplus \R_1(\u(p,q))\simeq\su(p,q)\oplus\Real J=\u(p,q).$$
Here by $V_\Lambda$ we denote the irreducible representations of
the Lie algebra $\sl(n,\Co)$ with the highest weight $\Lambda$. To
describe the above isomorphism, note that if $A\in\u(p,q)$, then
$H=JA\in\odot^2\Real^{2p,2q}$ is a symmetric endomorphism
commuting with $J$, or, equivalently, $H$ defines a symmetric
bilinear form such that $$H(JX,Y)+H(X,JY)=0.$$ Any $H$ with such
property may be obtained in this way. Next, such $H$ defines
$R_H\in\R(\u(p,q))$ given by
$$R_H(X,Y)=HX\wJ Y+X\wJ HY+2g(HJX,Y)J+2g(JX,Y)JH.$$
We set the notation $$H\owJ g=R_H.$$ Now, $\R'(\u(p,q))$ is
spanned by elements $R_H$ with trace-free $H$, and $J\in\u(p,q)$
defines the element $-\Id\owJ g\in \R_1(\u(p,q))$. Note that
$$\left(\frac{\Id}{2}\owJ g\right)(X,Y)=X\wJ Y+2g(JX,Y)J.$$
The above decomposition of $\R(\u(p,q))$ shows that any
$R\in\R(\u(p,q))$ may be written in the form
$$R=B+H\owJ g+\lambda\frac{\Id}{2}\owJ g,$$
where $B$ is a totally trace-free tensor, $H$ is trace-free and
$\lambda\in\Real$. It can be computed that
$H=-\frac{1}{2(n+2)}\Ric^0$ and $\lambda=-\frac{s}{4n(n+1)}$, here
$\Ric^0=\Ric-\frac{s}{2n}g$ is the trace-free part of the Ricci
operator $\Ric$ defined by $R$, and $s$ is the trace of $\Ric$, i.e.
the scalar curvature defined by $R$.

Let $(M,g,J)$ be a pseudo-K\"ahlerian manifold. Then its curvature
tensor admits the decomposition: \be
R=B-\frac{1}{2(n+2)}\left(\Ric-\frac{s}{2n}\id\right)\owJ
g-\frac{s}{4n(n+1)}\frac{\id}{2}\owJ g.\ee The tensor $B$ is
called the Bochner curvature tensor, and it is an analog of the
Weyl conformal tensor $W$. The above equality can be rewritten in
the form \be\label{RBK} R=B+K\owJ g,\ee where
$$K=-\frac{1}{2(n+2)}\left(\Ric-\frac{s}{4(n+1)}g\right)$$
is the analog of the Schouten tensor $L$.

Define the tensor $D$ in the following way: \be
D(Z,X,Y)=g(g^{ab}\nabla_{X_a}R(Z,X_b)Y,X)-\frac{1}{4(n+1)}g\left(\left(\frac{\id}{2}\owJ
g\right)(\grad s, Z)Y,X\right). \ee The tensor $D$ is the analog
of the Cotton tensor $C$. The following holds \cite{Ta67}:
\be\label{DB}
D(Z,X,Y)=\frac{n+2}{n}g(g^{ab}\nabla_{X_a}B(Z,X_b)Y,X). \ee Let
\be\label{tDD} g(\tilde D(X)Y,Z)=D(X,Y,Z),\ee then
$$\tilde D(X)=-\frac{n+2}{n}g^{ab}\nabla_{X_a}B(Z,X_b) $$
or, in other terms, \be\label{oprD} \tilde
D(Z)=-g^{ab}\nabla_{X_a}R(Z,X_b)+\frac{1}{4(n+1)}\left(\frac{\id}{2}\owJ
g\right)(\grad s, Z).\ee

\section{Decomposition of the covariant derivative of the curvature
tensor of a pseudo-K\"ahlerian manifold}\label{secnRK}

Let $U$ be an irreducible $\u(p,q)$-module such that $J$ acts as a
complex structure on $U$. Then the $\gl(n,\Co)$-module $U\otimes
\Co$ decomposes into the direct sum $V\oplus\bar V$ of irreducible
$\gl(n,\Co)$-modules, where $V$ and $\bar V$ are the eigenspaces
of the extension of $J$ to $U\otimes\Co$ corresponding to the
eigenvalues $\rm i$ and $-\rm i$, respectively. The representation
of $\gl(n,\Co)$ on $V\oplus\bar V$ is given in the matrix form \be
\label{formrep}\left.\left\{\left(\begin{array}{cc}\rho(A)&0\\0&-\overline{
\rho(A)}^t\end{array}\right)\right|A\in\gl(n,\Co)\right\},\ee
where $\rho:\gl(n,\Co)\to \gl(V)$ is the representation of
$\gl(n,\Co)$ on $V$. Note that $U$ is isomorphic to $V$ as the
complex vector space. The vector space $V$ is the highest weight
$\gl(n,\Co)$-module with a highest weight $\Lambda$, and we denote
$V$ by $V_\Lambda$. Thus, $U\simeq V_\Lambda$. Next, if $n\geq 3$
and $V_\Lambda$ is an irreducible representation of $\gl(n,\Co)$
with the highest weight $\Lambda$ and the labels of $\Lambda$ are
situated not symmetrically on the Dynkin diagram of $\sl(n,\Co)$,
then the restriction of this representation to $\u(p,q)$ is
irreducible \cite{V-O}.

The space $\R^\nabla(\u(p,q))$ admits the natural complex
structure, $$(J\cdot S)_X=S_{-JX}$$ given by the representation of
$\u(p,q)$ on $\R^\nabla(\u(p,q))$. Hence, each irreducible
component of $\R^\nabla(\u(p,q))$ is of the form $V_\Lambda$.

\begin{theorem}\label{Th3} If $n\geq 3$, then the $\u(p,q)$-module
  $\R^\nabla(\u(p,q))$ admits the following decomposition into the sum of
 irreducible $\u(p,q)$-modules:
$$\R^\nabla(\u(p,q))=\Q_0\oplus\Q'\oplus
\Q_1\simeq V_{3\pi_1+2\pi_{n-1}}\oplus V_{\pi_2+\pi_{n-1}}\oplus
\Real^{2p,2q}.$$

If $n=2$, then
$$\R^\nabla(\u(p,q))=\Q_0\oplus\Q'\oplus
\Q_1\simeq \Co^6\oplus\Co^4\oplus\Real^{2p,2q}.$$
\end{theorem}

Note that $$\Q'\simeq \P_0(\u(p,q))\simeq
(\odot^2(\Co^n)^*\otimes\Co^n)_0,$$ where
$(\odot^2(\Co^n)^*\otimes\Co^n)_0$ is the subspace of
$\odot^2(\Co^n)^*\otimes\Co^n$ consisting of tensors such that the
contraction of the upper index with any down index gives zero.
 The explicit isomorphism will
be constructed in the proof of Theorem \ref{Th3}.

Recall that we consider a pseudo-K\"ahlerian manifold  $(M,g,J)$
of signature $(2p,2q)$,  $p+q=n\geq 2$. Now we give the explicit
form of the above decomposition for the covariant derivative
$\nabla R$  of the  curvature tensor $R$ of $(M,g,J)$. Let us
first set some notation. Define the following tensors:
\begin{align}
T_{X}=&\frac{1}{8(n+2)(n+1)}\Big( J\circ(X\wJ J\grad s)-
g(\grad s, X)\Id\Big),\\
(\psi_1)_X=&-\frac{1}{2(n+3)}J\tilde D(JX)\owJ g,\\
(\psi_2)_X(Y,Z)=&-\frac{1}{2(n+3)}\left(\tilde D\big((Y\wJ
Z)X\big)+X\wJ \big(\tilde D(Z)Y- \tilde D(Y)Z\big)\right).
\end{align}

\begin{theorem}\label{Th4} Let $(M,g,J)$ be a pseudo-K\"ahlerian manifold of signature $(2p,2q)$, $p+q=n\geq 2$.
Then the covariant derivative $\nabla R$ of the  curvature tensor
$R$ of $(M,g,J)$ admits the following decomposition:
\begin{equation}\label{decnablaRK} \nabla R=Q_0+Q'+Q_1,\end{equation}
where  \begin{align}
\label{Q0}Q_0&=\nabla B+\frac{1}{n+2}\psi_1+\psi_2,\\
\label{Q'}Q'&=\psi_1-\psi_2,\\
(Q_1)_X&=T_{X}\owJ g.\end{align} \end{theorem}

The fact that $Q_0$ satisfies the second Bianchi identity can be
expressed in the form \begin{multline} g\big((\nabla_X
B(Y,Z)+\nabla_Y B(Z,X)+\nabla_Z
B(X,Y))V,U\big)\\=\frac{1}{2(n+2)}\Big(D(U,X,Y)h(V,Z)+D(U,Y,Z)h(V,X)+D(U,Z,X)h(V,Y)
\\+D(h(Z,U)V,Y,X)+D(h(X,U)V,Z,Y)+D(h(Y,U)V,X,Z)
\\-2D(g(JX,Y)JZ+g(JY,Z)JX+g(JZ,X)JY,V,U)\Big).\end{multline}
This gives the second Bianchi identity for the Bochner tensor $B$,
similar to the  second Bianchi identity for the Weyl tensor
\eqref{2BW}. This identity written using the index notation is
proved directly  in \cite{Omachi03}.

Consider some traces for the obtained tensors. The tensor $Q_0$ is
totally trace-free. It holds
\begin{align}
 (Q')_{X_a}(X,Y,Z,X_b)g^{ab}&=D(Z,X,Y),\\
 (Q_1)_{X_a}(X,Y,Z,X_b)g^{ab}&=\frac{1}{4(n+1)}g\left(\left(\frac{\id}{2}\owJ
g\right)(\grad s,Z)Y,X\right),\\
 (Q')_{X}(X_a,Y,X_b,Z)g^{ab}&=-D(JX,JZ,Y),\\
 (Q_1)_{X}(X_a,Y,X_b,Z)g^{ab}&=\frac{1}{4(n+1)}g\left(\left(\frac{\id}{2}\owJ
g\right)(X,J\grad s)JY,Z\right).
\end{align}

This shows that
$$\nabla_X\Ric(Y,Z)=-D(JX,JZ,Y)+\frac{1}{4(n+1)}g\left(\left(\frac{\id}{2}\owJ
g\right)(X,J\grad s)JY,Z\right),$$ and the space of values of the
covariant derivatives of the Ricci tensor is decomposed into the
direct sum of two irreducible $\u(p,q)$-modules isomorphic,
respectively, to $(\odot^2(\Co^n)^*\otimes\Co^n)_0$ and
$\Real^{2p,2q}$.

We also get
$$\nabla_X\Ric(JY,Z)=D(JX,Y,Z)-\frac{1}{4(n+1)}g\left(\left(\frac{\id}{2}\owJ
g\right)(X,J\grad s)Y,Z\right),$$ this corresponds to the fact
that $\nabla\Ric(J\cdot,\cdot)$ satisfies
$$\nabla_X\Ric(JY,Z)+\nabla_Y\Ric(JZ,X)+\nabla_Z\Ric(JX,Y)=0,$$
i.e. $\nabla\Ric(J\cdot,\cdot)\in\P(\u(TM))$ and gives the
decomposition of $\nabla\Ric(J\cdot,\cdot)$ corresponding  to the
decomposition $\P(\u(p,q))=\P_0(\u(p,q))\oplus \P_1(\u(p,q))$ (see
below).

\section{Proof of Theorem \ref{Th1}}\label{secPTh1}
Consider a subalgebra $\h\subset\so(p,q)$ and its complexification
$\h\otimes\Co\subset\so(\Real^{p,q}\otimes\Co)$.
 We will use the facts that
\be\label{compl}\R(\h\otimes\Co)=\R(\h)\otimes\Co,\qquad
\R^\nabla(\h\otimes\Co)= \R^\nabla(\h)\otimes\Co.\ee Let
$V_\Lambda$ denote the irreducible $\so(n,\Co)$-module of
$\so(n,\Co)$ with the highest  weight $\Lambda$. If the
restriction of this representation of $\so(n,\Co)$ to
$\so(p,q)\subset \so(n,\Co)$ is reducible (this happens  if and
only if $n\geq 4m+2$ or $n=4m+2$ and the labels of $\Lambda$ on
the Dynkin diagram of $\so(4m+2,\Co)$ are situated symmetrically),
then as $\so(p,q)$-module $V_\Lambda$ decomposes in the direct sum
$V^{\so(p,q)}_\Lambda\oplus i V^{\so(p,q)}_\Lambda$ for some
irreducible $\so(p,q)$-module $V^{\so(p,q)}_\Lambda$. When it does
not lead to confuse, we denote $V^{\so(p,q)}_\Lambda$ simply by
$V_\Lambda$ (this denotation is used in Sections \ref{secR} and
\ref{secnR}).

First suppose that $n\geq 5$. Using \eqref{decRso}, we get
\begin{multline} \R^\nabla(\so(p,q))\subset\Real^{p,q}\otimes\R(\so(p,q))=\Real^{p,q}\otimes(V_{2\bar\pi_2}\oplus V_{2\pi_1}\oplus \Real)\\=
\big(V_{\pi_1+2\bar\pi_2}\oplus V_{\lambda}\oplus
V_{\pi_1+\bar\pi_2}\big)\oplus\big(V_{3\pi_1} \oplus
V_{\pi_1+\bar\pi_2}\oplus\Real^{p,q}\big)\oplus\Real^{p,q},
\end{multline} where $\lambda=4\pi_2$ if $n=5$, and $\lambda=\pi_2+\pi_3$ if $n\geq 6$.
  Consider the $\so(n,\Co)$-submodule
$V_{\pi_1+2\bar\pi_2}\subset\Co^n\otimes V_{2\bar\pi_2}$. This
submodule is the highest one. Let $v_{-m},...,v_{-1},v_1,...,v_m$
and $v_{-m},...,v_{-1},v_0,v_1,...,v_m$ be the standard bases of
$\Co^{2m}$ and $\Co^{2m+1}$, respectively. A highest weight vector
of the $\so(n,\Co)$-module
 $V_{2\bar\pi_2}$ has the form $(v_1\wedge v_2)\odot(v_1\wedge v_2)$ \cite{Al}. Hence,
 $$\xi=v_1\otimes((v_1\wedge v_2)\odot(v_1\wedge v_2))$$ is a
highest weight vector  of the module $V_{\pi_1+2\bar\pi_2}$.
Clearly, $\xi\in\R^\nabla(\so(n,\Co))$, consequently,
$V_{\pi_1+2\bar\pi_2}\subset\R^\nabla(\so(n,\Co))$. Suppose that
$n\geq 6$. Consider the submodule
$$V_{\pi_2+\pi_3}\subset\Co^n\otimes
V_{2\pi_2}=V_{\pi_1+2\pi_2}\oplus V_{\pi_2+\pi_3}\oplus
V_{\pi_1+\pi_2}.$$
 Consider the element
$$\xi=v_3\otimes((v_1\wedge v_2)\odot(v_1\wedge v_2))\in
\Co^n\otimes V_{2\pi_2},$$ which has weight $\pi_2+\pi_3$. It is
easy to check that $\xi\not\in\R^\nabla(\so(n,\Co))$. Note that
the module $V_{\pi_1+\pi_2}$ does not contain the weight space of
weight $\pi_2+\pi_3$. Hence, there exist $\xi_1\in
V_{\pi_1+2\pi_2}$ and $\xi_2\in V_{\pi_2+\pi_3}$ such that
$\xi=\xi_1+\xi_2$. We obtain that $\xi_2\neq 0$ and
$\xi_2\not\in\R^\nabla(\so(n,\Co))$ this shows that
$V_{\pi_2+\pi_3}\not\subset\R^\nabla(\so(n,\Co))$. Similarly, if
$n=5$, then $V_{4\pi_2}\not\subset\R^\nabla(\so(n,\Co))$.

For any linear map $S:\Real^{p,q}\to\R(\so(p,q))$ define the map
\begin{equation} P:\Real^{p,q} \to\so(p,q),\quad
P(X)=\tr_{(1,5)}S_\cdot(\cdot,\cdot,X,\cdot)= g
^{ij}S_{X_i}(X,X_j).\end{equation}
 It is easy to check that for  any $S\in
\R^\nabla(\so(p,q))$, the element $\tr_{(1,5)}S$ belongs to
$\P(\so(p,q))$. For any subalgebra $\h\subset\so(p,q)$ and
$P\in\P(\h)$ define the vector
\begin{equation}\label{tRic}\tRic(P)= g ^{ij}P(X_i)X_j\in\Real^{p,q}.\end{equation}
In \cite{onecomp} it is shown that  the space $\P(\h)$ admits the
decomposition
$$\P(\h)=\P_0(\h)\oplus \P_1(\h),$$
where $\P_0(\h)=\ker\tRic$ and $\P_1(\h)$ is its orthogonal
complement.

It holds $$\P_0(\so(p,q))\simeq V_{\pi_1+\bar\pi_2},\quad
\P_1(\so(p,q))\simeq \Real^{p,q}.$$ In correspondence with this
decomposition, any   tensor $P\in \P(\so(p,q))$ can be decomposed
in the sum
\begin{equation}P=P_0+P_1,\quad
P_0(X)=P(X)+\frac{1}{n-1}\tRic(P)\wedge X,\quad
P_1(X)=-\frac{1}{n-1}\tRic(P)\wedge X.\end{equation} For any
tensor $P$ define the tensors $\Phi_1(P)$ and $\Phi_2(P)$   as
follows:
\begin{equation}\label{phi1}
\Phi_1(P)_X(Y,Z)=H_XY\wedge Z + Y\wedge H_XZ,\end{equation} where
$H:\Real^{p,q}\to\odot^2\Real^{p,q}$ is given by
$$ g (H_XY,Z)= g (P(Y)Z+P(Z)Y,X),$$
 and \begin{equation}\label{phi2}
\Phi_2(P)_X(Y,Z)= P\big((Y\wedge Z)X\big)+X\wedge
(P(Z)Y-P(Y)Z).\end{equation} It can be checked that $\Phi_1(P)$
and $\Phi_2(P)$ map $\Real^{p,q}$ to $\R(\so(p,q)$, and the maps
$\Phi_1$ and $\Phi_2$ are $\so(p,q)$-equivariant. Let
$\alpha,\beta\in\Real$ and $P\neq 0$. It can be shown that
$(\alpha \Phi_1+\beta\Phi_2)(P)$ belongs to $\R^\nabla(\so(p,q))$
if and only if $\beta =3\alpha$. Next,
\begin{equation}\label{trPhi}\tr_{(1,5)}(\Phi_1(P_0))=-3 P_0,\quad
\tr_{(1,5)}(\Phi_1(P_1))=(n-4)P_1,\quad
\tr_{(1,5)}(\Phi_2(P))=-nP,\end{equation} where
$P\in\P(\so(p,q))$,  $P_0\in\P_0(\so(p,q))$ and
$P_1\in\P_1(\so(p,q))$. This shows that $\R^\nabla(\so(p,q))$
contains exactly one submodule isomorphic to $V_{\pi_1+\bar\pi_2}$
and exactly one submodule isomorphic to $\Real^{p,q}$.

Recall that $\R(\so(p,q))$ contains the submodule
$\odot^2\Real^{p,q}$, hence $\Real^{p,q}\otimes \R(\so(p,q))$
contains the submodule $\odot^3\Real^{p,q}\simeq V_{3\pi_1}\oplus
\Real^{p,q}$. The inclusion $$ \odot^3\Real^{p,q}\hookrightarrow
\Real^{p,q}\otimes \R(\so(p,q))$$ is given by
$$T\mapsto S,\quad S_X=T_X\owedge  g .$$
Clearly, $S\in\R^\nabla(\so(p,q))$.  This shows that
$V_{3\pi_1}\subset \R^\nabla(\so(p,q))$ and it gives the explicit
inclusion $ \odot^3\Real^{p,q}\hookrightarrow
\R^\nabla(\so(p,q))$. For $n=3 $ and $4$ the proof is similar. We
have proved Theorem \ref{Th1}. $\Box$

\section{Proof of Theorem \ref{Th2}}\label{secPTh2}

Let $(M,g)$ be a pseudo-Riemannian manifold. Suppose that $n\geq
6$.  From Theorem \ref{Th1} it follows that
$$\nabla R=S_0'+S_0''+S'+S_1$$  for some tensors $S_0',S_0'',S',S_1$ such that
the values of these tensors at any $x\in M$ belong to the
submodules of $\R^\nabla(\so(T_xM))$ isomorphic respectively to
$V_{\pi_1+2\pi_2}$, $V_{3\pi_1}$, $V_{\pi_1+\pi_2}$ and
$\Real^{p,q}$.

From \eqref{decR1} it follows that  $$\nabla_X R=\nabla_X
W-\nabla_XL\owedge g.$$ Let $x\in M$. Then for  each $X\in T_xM$,
$\nabla_X W_x$ belongs to $V_{2\pi_2}\subset \R(\so(p,q))$, while
$\nabla_X L_x\owedge g_x$ belongs to
$\odot^2\Real^{p,q}=V_{2\pi_1}\oplus\Real\subset\R(\so(p,q))$. We
get
\begin{align} \nabla W_x&\in \Real^{p,q}\otimes V_{2\pi_2}=
V_{\pi_1+2\pi_2}\oplus V_{\pi_1+\pi_2}\oplus V_{\pi_2+\pi_3},\\
\label{razlnabL}\nabla L_x&\in  \Real^{p,q}\otimes
(V_{2\pi_1}\oplus\Real)= (V_{3\pi_1}\oplus
V_{\pi_1+\pi_2}\oplus\Real^{p,q})\oplus\Real^{p,q}.
\end{align}

Since $\R^\nabla(\so(p,q))$ does not contain the submodule
$V_{\pi_2+\pi_3}$, we obtain
\begin{equation}\label{razlnabW}\nabla W_x\in V_{\pi_1+2\pi_2}\oplus V_{\pi_1+\pi_2}.\end{equation}
From the above we get that there are unique elements
$P_{0x}\in\P_0(\so(p,q))$  and $P_{1x}\in\P_1(\so(p,q))$ such that
$$\tr_{(1,5)}(\nabla R_x-(\Phi_1+3\Phi_2)(P_{0x}+P_{1x}))=0.$$
Let us find these elements. Using the second Bianchi identity, it
can be shown that
$$g^{ij}\nabla_{X_i}
R(X,Y,Z,X_j)=C(Z,X,Y)+\frac{1}{2(n-1)}g(\grad s, (X\wedge Y)Z).$$
From this and \eqref{trPhi} we get
$$g(P_0(X)Y,Z)=\frac{1}{3(n+1)}C(X,Y,Z),\quad P_1(X)=-\frac{1}{4(n+2)(n-1)}\grad s\wedge X.$$
Let $P=P_0+P_1$. We obtain $$\nabla R_x-(\Phi_1+3\Phi_2)(P_x)\in
V_{\pi_1+2\pi_2}\oplus V_{3\pi_1}.$$ This means that
$$S_1=(\Phi_1+3\Phi_2)(P_1),\quad S'=(\Phi_1+3\Phi_2)(P_0), \quad
S'_0+S''_0= \nabla R-(\Phi_1+3\Phi_2)(P).$$ Note that
$\Phi_1(P_0)=\varphi_1$ and $\Phi_2(P_0)=\varphi_2$. Recall that
$\odot^3\Real^{p,q}\simeq V_{3\pi_1}\oplus\Real^{p,q}$. This,
Theorem \ref{Th1} and \eqref{razlnabL} show that
$$(S''_0+S_1)_X=T_X\owedge g,$$ where
$$g(T_XY,Z)=-\frac{1}{3}\left(\nabla_XL(Y,Z)+\nabla_YL(Z,X)+\nabla_ZL(X,Y)\right)$$
is the symmetrization of $-\nabla L$. Using this and the
expression for $S_1$, we find $S''_0$. Finally we find $S'_0$ from
the equality
$$S'_0=\nabla W-\nabla L\owedge g-S''_0-S_1-S'.$$

 Theorem \ref{Th2} is proved. $\Box$

\section{Proof of Theorem \ref{Th3}}\label{secPTh3}

First consider the complexification of the space $\R(\u(p,q))$. As
it is noted above, the representation of
$\gl(n,\Co)=\u(p,q)\otimes\Co$ on the space
$\Real^{2p,2q}\otimes\Co$ decomposes into the sum $V\oplus\bar
V=V_{\pi_1}\oplus V_{\pi_{n-1}}$ and the representation is of the
form \eqref{formrep}. From this, \eqref{compl} and the Bianchi
identity it follows that any $R\in\R(\u(p,q))\otimes\Co$ satisfies
$R(V,V)=R(\bar V,\bar V)=0$. Consequently,
$$\R(\u(p,q))\otimes\Co\subset (V_{\pi_1}\otimes
V_{\pi_{n-1}})\otimes\gl(n,\Co).$$ More precisely,
$$\R(\u(p,q))\otimes\Co=V_{2\pi_1+2\pi_{n-1}}\oplus\sl(n,\Co)\oplus\Co.$$
Hence,
\begin{multline} \R^\nabla(\u(p,q))\otimes\Co\subset (V_{\pi_1}\oplus V_{\pi_{n-1}})\otimes(\R(\u(p,q))\otimes\Co)
\\=(V_{3\pi_1+2\pi_{n-1}}\oplus V_{\pi_1+\pi_{2}+2\pi_{n-1}}\oplus
2V_{2\pi_1+\pi_{n-1}}\oplus 2V_{\pi_1}\oplus
V_{\pi_2+\pi_{n-1}})\\\oplus (V_{2\pi_1+3\pi_{n-1}}\oplus
V_{2\pi_1+\pi_{n-2}+\pi_{n-1}}\oplus 2V_{\pi_1+2\pi_{n-1}}\oplus
2V_{\pi_{n-1}}\oplus V_{\pi_1+\pi_{n-2}}).
\end{multline}
Since $\R^\nabla(\u(p,q))$ admits a complex structure, each
irreducible submodule of $\R^\nabla(\u(p,q))$ defines two
$\gl(n,\Co)$-irreducible submodules of $\R^\nabla(\u(p,q))\otimes
\Co$.

  Consider the $\gl(n,\Co)$-submodule
$V_{3\pi_1+2\pi_{n-1}}\subset\Co^n\otimes V_{2\pi_1+2\pi_{n-1}}$.
This submodule is the highest one. Let $v_1,...,v_n$ and
$v_{-1},...,v_{-n}$ be the basis of $V=V_{\pi_1}=\Co^{n}$ and the
dual basis of $\bar V=V_{\pi_{n-1}}=(\Co^{n})^*$, respectively. A
highest weight vector of the $\gl(n,\Co)$-module
 $V_{2\pi_1+2\pi_{n-1}}$ has the form $(v_1\otimes v_{-n})\otimes (v_1\otimes v_{-n})$. Hence,
 $$\xi=v_1\otimes(v_1\otimes v_{-n})\otimes (v_1\otimes v_{-n})$$ is a
highest weight vector  of the module $V_{3\pi_1+2\pi_{n-1}}$.
Clearly, $\xi\in\R^\nabla(\u(p,q))\otimes\Co$, consequently,
$V_{3\pi_1+2\pi_{n-1}}\subset\R^\nabla(\u(p,q))\otimes\Co$.
Similarly,
$V_{2\pi_1+3\pi_{n-1}}\subset\R^\nabla(\u(p,q))\otimes\Co$. These
two modules define an irreducible submodule of
$\R^\nabla(\u(p,q))$ that we denote by $\Q_0$.  Consider the
submodule
$$V_{\pi_1+\pi_2+2\pi_{n-1}}\subset\Co^n\otimes
V_{2\pi_1+2\pi_{n-1}}=V_{3\pi_1+2\pi_{n-1}}\oplus
V_{\pi_1+\pi_{2}+2\pi_{n-1}}\oplus V_{2\pi_1+\pi_{n-1}}.$$
 Consider the element
$$\xi=v_2\otimes(v_1\otimes v_{-n})\otimes(v_1\otimes v_{-n})\in
\Co^n\otimes V_{2\pi_1+2\pi_{n-1}},$$ which has weight
$\pi_1+\pi_2+2\pi_{n-1}$. It is easy to check that
$\xi\not\in\R^\nabla(\u(p,q))\otimes\Co$. Note that the module
$V_{2\pi_1+\pi_{n-1}}$ does not contain the weight space of the
weight $\pi_1+\pi_2+2\pi_{n-1}$. Hence, there exist $\xi_1\in
V_{3\pi_1+2\pi_{n-1}}$ and $\xi_2\in V_{\pi_1+\pi_2+2\pi_{n-1}}$
such that $\xi=\xi_1+\xi_2$. We obtain that $\xi_2\neq 0$ and
$\xi_2\not\in\R^\nabla(\u(p,q))\otimes\Co$ this shows that
$V_{\pi_1+\pi_2+2\pi_{n-1}}\not\subset\R^\nabla(\u(p,q))\otimes\Co$.
Similarly,
$V_{2\pi_1+\pi_{n-2}+\pi_{n-1}}\not\subset\R^\nabla(\u(p,q))\otimes\Co$.

The submodule $$2V_{2\pi_1+\pi_{n-1}}\oplus
2V_{\pi_1+2\pi_{n-1}}\subset
(\Real^{2p,2q}\otimes\Co)\otimes(\R(\u(p,q))\otimes\Co)$$ defines
two irreducible submodules in $\Real^{2p,2q}\otimes \R(\u(p,q))$
isomorphic to the irreducible $\u(p,q)$-module
$V_{\pi_1+2\pi_{n-1}}$. Similarly, the submodule
$$2V_{\pi_1}\oplus 2V_{\pi_{n-1}}\subset
(\Real^{2p,2q}\otimes\Co)\otimes(\R(\u(p,q))\otimes\Co)$$ defines
two irreducible submodules in $\Real^{2p,2q}\otimes \R(\u(p,q))$
isomorphic to $\Real^{2p,2q}$. We will show that
$\R^\nabla(\u(p,q))$ contains only one copy of each of these
submodules, we denote them by $\Q'$ and $\Q_1$, respectively. For
this we turn to the space $\P(\u(p,q))$. From the results of
Leistner \cite{Leistner}, it follows that
$$\P(\u(p,q))\otimes\Co\simeq (\gl(n,\Co)\subset V)^{(1)}\oplus
(\gl(n,\Co)\subset \bar V)^{(1)},$$ where $(\gl(n,\Co)\subset
V)^{(1)}$ denotes the first prolongation. It holds,
$$(\gl(n,\Co)\subset V)^{(1)}\simeq
\odot^2(\Co^n)^*\otimes\Co^n\simeq V_{\pi_1+2\pi_{n-1}}\oplus
V_{\pi_{n-1}}.$$ Thus, $$\P(\u(p,q))\simeq
\odot^2(\Co^n)^*\otimes\Co^n.$$

This  isomorphism has the following explicate form given in
\cite{onecomp}. Let $$S\in
\odot^2(\Co^n)^*\otimes\Co^n\subset(\Co^n)^*\otimes \gl(n,\Co).$$
We fix an identification $$\Co^n=\Real^{2p,2q}=\Real^{p,q}\oplus
i\Real^{p,q}$$ and choose an orthonormal basis $e_1,...,e_n$ of
$\Real^{p,q}$. Define the complex numbers $S_{abc}$,
$a,b,c=1,...,n$ such that $$S(e_a)e_b=\sum_c S_{acb}e_c.$$ It
holds $S_{abc}=S_{cba}$. Define a map
$S_1:\Real^{2p,2q}\to\gl(2n,\Real)$ by the conditions
$$S_1(e_a)e_b=\sum_c \overline {S_{abc}}e_c,\quad
S_1(ie_a)=-iS_1(e_a),\quad S_1(e_a)ie_b=iS_1(e_a)e_b.$$ The map
$P=S-S_1:\Real^{2n}\to\gl(2n,\Real)$ belongs to $\P(\u(p,q))$ and
any element of $\P(\u(p,q))$ is of this form. Such element belongs
to $\P_0(\u(p,q))$ if and only if $\sum_b S_{abb}=0$ for all
$a=1,...,n$, i.e. $S\in (\odot^2(\Co^n)^*\otimes\Co^n)_0$. Next,
$$\P(\u(p,q))=\P_0(\u(p,q))\oplus\P_1(\u(p,q))\simeq
(\odot^2(\Co^n)^*\otimes\Co^n)_0\oplus\Real^{p,q}.$$ According to
this, any $P\in\P(\u(p,q))$ admits the decomposition \be
P=\left(P-\frac{1}{2(n+1)}\left(\frac{\id}{2}\owJ
g\right)(\tRic(P),\cdot)\right)+\frac{1}{2(n+1)}\left(\frac{\id}{2}\owJ
g\right)(\tRic(P),\cdot).\ee Any linear map
$Q:\Real^{2p,2q}\to\R(\u(p,q))$ defines the map
\begin{equation} P:\Real^{2p,2q} \to\u(p,q),\quad
P(X)=(\tr_{(1,5)})Q_\cdot(\cdot,\cdot,X,\cdot)= g
^{ab}Q_{X_a}(X,X_b).\end{equation}
 If $Q\in
\R^\nabla(\u(p,q))$, then the element $\tr_{(1,5)}Q$ belongs to
$\P(\u(p,q))$. For any $P\in\P(\u(p,q))$  define the tensors
$\Psi_1(P)$ and $\Psi_2(P)$   as follows:
\begin{equation}\label{psi1}
\Psi_1(P)_X=JP(JX)\owJ g,\end{equation}
\begin{equation}\label{psi2} \Psi_2(P)_X(Y,Z)= P\big((Y\wJ
Z)X\big)+X\wJ (P(Z)Y-P(Y)Z).\end{equation} It can be checked that
$\Psi_1(P)$ and $\Psi_2(P)$ map $\Real^{2p,2q}$ to $\R(\u(p,q))$,
and the maps $\Psi_1$ and $\Psi_2$ are $\u(p,q)$-equivariant. It
can be shown that if $P\neq 0$, then $\Psi_1(P)\not\in
\R^\nabla(\u(p,q))$. Moreover,
$(\Psi_1-\Psi_2)(P)\in\R^\nabla(\u(p,q))$ for any
$P\in\P(\u(p,q))$.
 Next,
\begin{equation}\label{trPsi}\tr_{(1,5)}((\Psi_1-\Psi_2)P_0)=2(n+3) P_0,\quad
\tr_{(1,5)}((\Psi_1-\Psi_2)P_1)=4(n+2) P_1,\end{equation} where
$P_0\in\P_0(\u(p,q))$ and $P_1\in\P_1(\u(p,q))$. This shows that
$\R^\nabla(\u(p,q))$ contains exactly one submodule isomorphic to
$(\odot^2(\Co^n)^*\otimes\Co^n)_0$ and exactly one submodule
isomorphic to $\Real^{2p,2q}$.

We are left with the submodules $V_{\pi_2+\pi_{n-1}}$ and
$V_{\pi_1+\pi_{n-2}}$ of
$(\Real^{2p,2q}\otimes\Co)\otimes(\R(\u(p,q))\otimes\Co)$. These
submodules define an irreducible submodule in
$\Real^{2p,2q}\otimes\R(\u(p,q))$. This submodule is contained in
$\Real^{2p,2q}\otimes\R'(\u(p,q))$. Hence any element of this
submodule is of the form $$\xi_X=\Upsilon_X\owJ g,$$ where for
each $X$, $\Upsilon_X$ is symmetric and commutes with $J$.
Clearly, the contraction of any two indexes of $\Upsilon$ gives
$0$. Suppose that $\xi\in\R^\nabla(\u(p,q))$. Then,
$$g^{ab}g((\xi_X(X_a,Z)+\xi_{X_a}(Z,X)+\xi_Z(X,X_a))X_b,V)=0.$$
This implies
$$(2n+3)\Upsilon_X(Z,V)-(2n+3)\Upsilon_Z(X,V)-\Upsilon_{JX}(JZ,V)+\Upsilon_{JZ}(JX,V)-2\Upsilon_{JV}(JZ,X)=0.$$
Simple manipulations show
$$\Upsilon_V(X,Z)=\Upsilon_{JV}(X,JZ).$$ Finally,
$$\Upsilon_V(X,Z)=\Upsilon_{JV}(X,JZ)=-\Upsilon_{JV}(JX,Z)=-\Upsilon_{V}(X,Z),$$
and $\Upsilon=0$. Thus the submodule under the consideration is
not contained in $\R^\nabla(\u(p,q))$. The theorem is proved.
$\Box$

\section{Proof of Theorem \ref{Th4}}\label{secPTh4}

Let $(M,g,J)$ be a pseudo-K\"ahlerian manifold.  From Theorem
\ref{Th3} it follows that
$$\nabla R=Q_0+Q'+Q_1$$  for some tensors $Q_0,Q',Q_1$ such that
the values of these tensors at any $x\in M$ belong to the
submodules of $\R^\nabla(\u(T_xM))$ isomorphic respectively to
$V_{3\pi_1+2\pi_{n-1}}$, $V_{\pi_1+2\pi_{n-1}}$ and $\Real^{p,q}$.

From \eqref{RBK} it follows that  $$\nabla_X R=\nabla_X
B+\nabla_XK\owJ g.$$ Let $x\in M$. Then for  each $X\in T_xM$,
$\nabla_X B_x$ belongs to $V_{2\pi_1+2\pi_{n-1}}\subset
\R(\u(p,q))$, while $\nabla_X K_x\owJ g_x$ belongs to
$\u(p,q)\subset\R(\u(p,q))$. We get
\begin{align} \nabla B_x&\in \Real^{2p,2q}\otimes V_{2\pi_1+2\pi_{n-1}}=
V_{3\pi_1+2\pi_{n-1}}\oplus V_{\pi_1+\pi_2+2\pi_{n-1}}\oplus V_{\pi_1+2\pi_{n-1}},\\
\label{razlnabK}\nabla K_x&\in  \Real^{2p,2q}\otimes\u(p,q) =
(V_{\pi_1+2\pi_{n-1}}\oplus
V_{\pi_2+\pi_{n-1}}\oplus\Real^{2p,2q})\oplus\Real^{2p,2q}.
\end{align}

Since $\R^\nabla(\u(p,q))$ does not contain the submodule
$V_{\pi_1+\pi_2+2\pi_{n-1}}$, we obtain
\begin{equation}\label{razlnabB}\nabla B_x\in V_{3\pi_1+2\pi_{n-1}}\oplus  V_{\pi_1+2\pi_{n-1}}.\end{equation}
 Since $$g^{ab}\nabla_{X_a}
R_x(\cdot,X_b)\in\P(\u(p,q)),$$ this element can be decomposed in
the form \be\label{lab64} g^{ab}\nabla_{X_a}
R_x(\cdot,X_b)=-\tilde D_x+\left(\frac{\id}{2}\owJ
g\right)(V,\cdot)\ee for some $\tilde D_x\in\P_0(\u(p,q))$ and
$V\in\Real^{2p,2q}$. From the above it follows that there are
unique elements $P_{0x}\in\P_0(\u(p,q))$ and
$P_{1x}\in\P_1(\u(p,q))$ such that
$$Q_{1x}=(\Psi_1-\Psi_2)P_{0x},\quad Q'_x=(\Psi_1-\Psi_2)P_{1x}.$$
Let us find these elements. It holds
$$\tr_{(1,5)}(\nabla R_x)=\tr_{(1,5)}Q_{1x}+\tr_{(1,5)}Q'_x=2(n+3)P_{0x}+4(n+2)P_{1x}.$$
We conclude that \be P_{0x}=-\frac{1}{2(n+3)}\tilde D_x,\quad
P_{1x}=\frac{1}{4(n+2)}\left(\frac{\id}{2}\owJ
g\right)(V,\cdot).\ee The last equality implies
$$\tRic(P_{1x})=-\frac{n+1}{2(n+2)}V.$$
Taking a trace in \eqref{lab64}, we get $$\grad s=4(n+1)V.$$
Summarizing all the above, we obtain that $$
P_{1}=\frac{1}{16(n+2)(n+1)}\left(\frac{\id}{2}\owJ g\right)(\grad
s,\cdot),$$ and $\tilde D$ satisfies \eqref{tDD}, \eqref{oprD},
\eqref{DB}. Thus we have found $Q_1$ and $Q'$. The component $Q_0$
can be found as $$Q_0=\nabla R-Q_1-Q'.$$
 Theorem \ref{Th4} is proved. $\Box$

\section{Pseudo-Riemannian manifolds with other irreducible
holonomy algebras}\label{secOther}

In the same way as above, the spaces  $\R^\nabla(\h)$ may be
computed also for non-locally symmetric pseudo-Riemannian
manifolds with irreducible holonomy algebras $\h$ different from
$\so(p,q)$ and $\u(p,q)$. Most of these manifolds are Ricci-flat
or Einstein \cite{Bryant2}, consequently $\nabla R$ coincides with
$\nabla W$. In the following theorem the results of the
computations are given. The proof of the theorem is similar to the
proofs of Theorems \ref{Th1} and \ref{Th3}, by this reason we omit
it.

\begin{theorem} It holds
\begin{align*}
\R^\nabla(\su(p,q))&=\Q_0\simeq V^\Co_{3\pi_1+2\pi_{n-1}}, \text{ if } p+q\geq 2,\quad \R^\nabla(\su(1,1))\simeq V^\Co_{6\pi_1}\\
\R^\nabla(\sp(p,q)\oplus\sp(1))\otimes\Co&=\R^\nabla(\sp(p,q))\otimes\Co\simeq
V^\Co_{5\pi_1}, \text{ if } p+q\geq 2\\
\R^\nabla(\sp(2m,\Real)\otimes\sl(2,\Real))\otimes \Co&\simeq
V_{5\pi_1}, \text{ if } m\geq 2\\
\R^\nabla(\h)\otimes\Co&\simeq V_{\pi_1+2\pi_2},\text{ if }
\h=\g_2\subset\so(7),\,\g_{2(2)}\subset\so(3,4),\\
\R^\nabla(\h)\otimes\Co&\simeq V_{2\pi_2+\pi_3},\text{ if }
\h=\spin(7)\subset\so(8),\,\spin(3,4)\subset\so(4,4).
\end{align*}
If $\h\subset\so(n,n)$ is one of the holonomy algebras
$\so(n,\Co)\subset\so(n,n)$,
$\sp(2m,\Co)\oplus\sl(2,\Co)\subset\so(4m,4m)$,
$\g_2^\Co\subset\so(7,7)$, $\spin(7,\Co)\subset\so(7,7)$, then
$\R^\nabla(\h)$ is isomorphic to $\R^\nabla(\h_0)\otimes\Co$,
where $\h_0\subset\so(n)$ is one of $\so(n)$,
$\sp(m)\oplus\sp(1)$, $\g_2$, $\spin(7)$, respectively.

\end{theorem}

\end{document}